\documentclass[twocolumn]{autart}    

\usepackage{graphicx}          

\usepackage{amsmath}
\usepackage{amssymb}

\def\R{\mathbb{R}}

\usepackage{graphics} 
\usepackage{epsfig} 
\usepackage{epstopdf}
\usepackage{subfigure}
\usepackage{color}

\graphicspath{{./Figures/}}

\allowdisplaybreaks[4]

\newcounter{cst}

\begin{document}

\begin{frontmatter}

\title{Predictor-Based Output Feedback Stabilization of an Input Delayed Parabolic PDE with Boundary Measurement\thanksref{footnoteinfo}} 

\thanks[footnoteinfo]{Corresponding author H.~Lhachemi. The work of C. Prieur has been partially supported by MIAI@Grenoble Alpes (ANR-19-P3IA-0003)
}

\author[CS]{Hugo Lhachemi}\ead{hugo.lhachemi@centralesupelec.fr},
\author[GIPSA-lab]{Christophe Prieur}\ead{christophe.prieur@gipsa-lab.fr}, 

\address[CS]{Universit{\'e} Paris-Saclay, CNRS, CentraleSup{\'e}lec, Laboratoire des signaux et syst{\`e}mes, 91190, Gif-sur-Yvette, France}  
\address[GIPSA-lab]{Universit{\'e} Grenoble Alpes, CNRS, Grenoble-INP, GIPSA-lab, F-38000, Grenoble, France}             

\begin{keyword}                           
Input delayed reaction-diffusion PDEs, predictor, output feedback, boundary control             
\end{keyword}                             

\begin{abstract}                          
This paper is concerned with the output feedback boundary stabilization of general 1-D reaction diffusion PDEs in the presence of an arbitrarily large input delay. We consider the cases of Dirichlet/Neumann/Robin boundary conditions for the both boundary control and boundary condition. The boundary measurement takes the form of a either Dirichlet or Neumann trace. The adopted control strategy is composed of a finite-dimensional observer estimating the first modes of the PDE coupled with a predictor to compensate the input delay. In this context, we show for any arbitrary value of the input delay that the control strategy achieves the exponential stabilization of the closed-loop system, for system trajectories evaluated in $H^1$ norm (also in $L^2$ norm in the case of a Dirichlet boundary measurement), provided the dimension of the observer is selected large enough. The reported proof of this result requires to perform both control design and stability analysis using simultaneously the (non-homogeneous) original version of the PDE and one of its equivalent homogeneous representations.
\end{abstract}

\end{frontmatter}

\section{Introduction}\label{sec: Introduction}

Since time delays are ubiquitous in practical applications, feedback control of finite-dimensional systems in the presence of input delays has been extensively studied~\cite{artstein1982linear,richard2003time}. The extension of this topic to Partial Differential Equations (PDEs) has attracted much attention in the recent years~\cite{hashimoto2016stabilization,nicaise2009stability,wang2018delay}.

This paper is concerned with the feedback stabilization of reaction-diffusion PDEs in the presence of an arbitrarily long input delay. One of the very first contributions on this topic was reported in~\cite{krstic2009control} using a backstepping control design technique (see also~\cite{krstic2009compensating} for the related problem of sensor dynamics governed by diffusion PDEs). More recently, the possibility to combine classical spectral reduction methods~\cite{coron2004global,coron2006global,russell1978controllability} (which are based on the fact that the associated eigenfunctions form a Riesz basis) and the design of a classical predictor feedback~\cite{artstein1982linear,bekiaris2012compensation,karafyllis2017predictor} on a finite-dimensional truncated model of the original PDE was reported in~\cite{prieur2018feedback} in the case of a state-feedback. Extensions of this approach in various directions, also in the context of state-feedback, were reported in~\cite{lhachemi2020feedback,lhachemi2021indomain,lhachemi2020pi}. One of the main advantages of spectral reduction methods for parabolic PDEs is that they allow the design of a finite-dimensional state-feedback, making them particularly relevant for practical applications. However, sole state-feedback control of PDEs is generally inapplicable in practice because the distributed nature of the state makes it essentially impossible to measure. Hence the design of an observer is generally required. Since the plant is a PDE, the observer itself generally takes the form of a PDE synthesized using a backstepping procedure~\cite{krstic2008boundary}. Hence the partial state-feedback can be coupled with the observer to ensure the stability of the closed-loop plant; see e.g.~\cite{katz2020boundary} where sufficient LMI conditions are derived with robustness aspects w.r.t. small enough delays. In order to avoid the pitfall of late lumping approximations required for the implementation of observers with infinite dimensional dynamics, a number of works have been devoted to the design of finite-dimensional observer-based control strategies for parabolic PDEs~\cite{balas1988finite,curtain1982finite,grune2021finite,harkort2011finite,katz2020constructive,lhachemi2020finitePI,lhachemi2020finite,lhachemi2021local,lhachemi2021nonlinear,sakawa1983feedback,sano2012stability}. In this work, we take advantage of the control architecture initially reported in~\cite{sakawa1983feedback} augmented with the LMI-based procedure introduced in~\cite{katz2020constructive}. More precisely, we leverage the enhanced procedure reported in~\cite{lhachemi2020finite} that extends for general reaction-diffusion PDEs the LMI-based approach reported in~\cite{katz2020constructive} to Dirichlet and/or Neumann boundary control and measurement (see also~\cite{katz2020finite} with a different approach but limited to Dirichlet measurements).

We address the finite-dimensional observer-based output feedback boundary stabilization of general 1-D reaction diffusion PDEs in the presence of an arbitrarily large input delay. A solution to this control design problem was reported for the first time in~\cite{katz2021sub} by combining a finite-dimensional observer and a predictor (used to compensate the input delay) in the very specific setting of a reaction-diffusion equation with Neumann boundary control, a bounded output operator, and with stability of the closed-loop system assessed in $L^2$ norm for arbitrarily large value of the input delay. However, the approach developed in~\cite{katz2021sub} and which solely relies on the original (non-homogeneous) representation of the PDE is strongly tailored for the above-mentioned setting and is hardly extendable to other types of boundary conditions (Dirichlet/Robin), to unbounded measurement operators (Dirichlet/Neumann), and to system trajectories evaluated in $H^1$ norm. This is because, introducing $(\beta_n)_{n \geq 1}$ the coefficients of projection of the boundary control operator into the Hilbert basis formed by the eigenstructures of the underlying Sturm-Liouville operator, the Neumann setting is such that $\beta_n = O(1)$ while the Dirichlet/Robin configurations give in general no better than $\beta_n = O(\sqrt{\lambda_n})$ where $\lambda_n$ are the eigenvalues of the problem that grow in $n^2$. This difference of asymptotic behavior has a major impact on the control design procedure since the proof reported in~\cite{katz2021sub}, which is limited to the case of a bounded measurement operator with trajectories evaluated in $L^2$ norm, heavily relies on the convergence of the series $\sum \beta_n^2 / \lambda_n$ that is only granted for Neumann boundary control. This issue becomes even more stringent in the case of an unbounded measurement operator and for system trajectories evaluated in $H^1$ norm, as it will be further detailed in Remark~\ref{rem: position of this work} of this paper after the introduction of the model of the system and its modal decomposition.

In this paper, we completely solve the control design problem of output feedback stabilization of general 1-D reaction-diffusion PDEs in the presence of an arbitrarily large input delay, with Dirichlet/Neumann/Robin boundary control/condition, and with a boundary measurement selected as a either Dirichlet or Neumann trace. The employed control architecture combines a finite-dimensional observer and a predictor. In the case of a Dirichlet (resp. Neumann) measurement, we assess the exponential stability of the closed-loop system in $L^2$ and $H^1$ norms (resp. $H^1$ norm) for arbitrarily large value of the input delay provided the dimension of the observer is selected large enough. This is achieved by leveraging an adequate scaling procedure~\cite{lhachemi2020finite} and by considering simultaneously both the original (non-homogeneous) representation of the PDE and one of its homogeneous versions obtain using a change of variable for control design and Lyapunov stability analysis. 

The paper is organized as follows. Notations and properties of Sturm-Liouville operators are presented in Section~\ref{sec: preliminaries}. The problem setting is presented in Section~\ref{sec: problem setting} followed by a preliminary spectral reduction of the problem. The case of a Dirichlet measurement is studied in Section~\ref{sec: Case of a Dirichlet measurement} while the case of a Neumann measurement is analyzed in Section~\ref{sec: Case of a Neumann measurement}. A numerical example is provided in Section~\ref{sec: numerical example}. Finally, concluding remarks are formulated in Section~\ref{sec: conclusion}.

\section{Notation and properties}\label{sec: preliminaries}

\subsection{Notation}

Spaces $\R^n$ are endowed with the Euclidean norm $\Vert\cdot\Vert$. The corresponding induced norms of matrices are also denoted by $\Vert\cdot\Vert$. For any two vectors $X$ and $Y$ of arbitrary dimensions, $ \mathrm{col} (X,Y)$ stands for the vector $[X^\top,Y^\top]^\top$. The space of square integrable functions on $(0,1)$ is denoted by $L^2(0,1)$ and is endowed with the usual inner product $\langle f , g \rangle = \int_0^1 f(x) g(x) \,\mathrm{d}x$ and with associated norm denoted by $\Vert \cdot \Vert_{L^2}$. For an integer $m \geq 1$, $H^m(0,1)$ denotes the $m$-order Sobolev space and is equipped with its usual norm $\Vert \cdot \Vert_{H^m}$. For a symmetric matrix $P \in\R^{n \times n}$, $P \succeq 0$ (resp. $P \succ 0$) means that $P$ is positive semi-definite (resp. positive definite).

\subsection{Properties of Sturm-Liouville operators}\label{subsec: prop sturm liouville operators}

Let $\theta_1,\theta_2\in[0,\pi/2]$, $p \in \mathcal{C}^1([0,1])$ and $q \in \mathcal{C}^0([0,1])$ with $p > 0$ and $q \geq 0$. Let the Sturm-Liouville operator $\mathcal{A} : D(\mathcal{A}) \subset L^2(0,1) \rightarrow L^2(0,1)$ be defined by $\mathcal{A}f = - (pf')' + q f$ on the domain $D(\mathcal{A}) = \{ f \in H^2(0,1) \,:\, c_{\theta_1} f(0) - s_{\theta_1} f'(0) = c_{\theta_2} f(1) + s_{\theta_2} f'(1) = 0 \}$ where $c_{\theta_i} = \cos\theta_i$ and $s_{\theta_i} = \sin\theta_i$. The eigenvalues $\lambda_n$, $n \geq 1$, of $\mathcal{A}$ are simple, non negative (due to $\theta_1,\theta_2\in[0,\pi/2]$ and $q \geq 0$), and form an increasing sequence with $\lambda_n \rightarrow + \infty$ as $n \rightarrow + \infty$. The corresponding unit eigenvectors $\phi_n \in L^2(0,1)$ form a Hilbert basis. The domain of the operator $\mathcal{A}$ is also characterized by $D(\mathcal{A}) = \{ f \in L^2(0,1) \,:\, \sum_{n\geq 1} \vert \lambda_n \vert ^2 \vert \left< f , \phi_n \right> \vert^2 < +\infty \}$. Let $p_*,p^*,q^* \in \R$ be such that $0 < p_* \leq p(x) \leq p^*$ and $0 \leq q(x) \leq q^*$ for all $x \in [0,1]$, then it holds 
$
0 \leq \pi^2 (n-1)^2 p_* \leq \lambda_n \leq \pi^2 n^2 p^* + q^*
$
for all $n \geq 1$~\cite{orlov2017general}. Moreover if $p \in \mathcal{C}^2([0,1])$, we have (see, e.g., \cite{orlov2017general}) that $\phi_n (\xi) = O(1)$ and $\phi_n' (\xi) = O(\sqrt{\lambda_n})$ as $n \rightarrow + \infty$ for any given $\xi \in [0,1]$. Assuming further that $q > 0$, an integration by parts and the continuous embedding $H^1(0,1) \subset L^\infty(0,1)$ show the existence of constants $C_1,C_2 > 0$ such that
\begin{align}
C_1 \Vert f \Vert_{H^1}^2 \leq 
\sum_{n \geq 1} \lambda_n \left< f , \phi_n \right>^2
= \left< \mathcal{A}f , f \right>
\leq C_2 \Vert f \Vert_{H^1}^2 \label{eq: inner product Af and f}
\end{align}
for any $f \in D(\mathcal{A})$. The latter inequalities and the Riesz-spectral nature of $\mathcal{A}$ imply that the series expansion $f = \sum_{n \geq 1} \left< f , \phi_n \right> \phi_n$ holds in $H^2(0,1)$ norm for any $f \in D(\mathcal{A})$. Due to the continuous embedding $H^1(0,1) \subset L^{\infty}(0,1)$, we obtain that $f(0) = \sum_{n \geq 1} \left< f , \phi_n \right> \phi_n(0)$ and $f'(0) = \sum_{n \geq 1} \left< f , \phi_n \right> \phi_n'(0)$. We finally define, for any integer $N \geq 1$, $\mathcal{R}_N f = \sum_{n \geq N+1} \left< f , \phi_n \right> \phi_n$.

\section{Problem setting and preliminary spectral reduction}\label{sec: problem setting}

\subsection{Problem setting}

We consider in this paper the input delayed reaction-diffusion system described by
\begin{subequations}\label{eq: PDE}
\begin{align}
& z_t(t,x) = (p(x) z_x(t,x))_x - \tilde{q}(x) z(t,x) \\
& c_{\theta_1} z(t,0) - s_{\theta_1} z_x(t,0) = 0 \\
& c_{\theta_2} z(t,1) + s_{\theta_2} z_x(t,1) = u(t-h) \\
& z(0,x) = z_0(x)
\end{align}
\end{subequations}
for $t > 0$ and $x \in (0,1)$ where $\theta_1,\theta_2 \in [0,\pi/2]$, $p \in\mathcal{C}^2([0,1])$ with $p > 0$, $\tilde{q} \in\mathcal{C}^0([0,1])$, and the input delay $h > 0$. Here $z(t,\cdot)$ is the state of the PDE at time $t$, $u(t-h)$ is the delayed version of the command $u(t)$, and $z_0$ is the initial condition. We assume throughout the paper that $u(\tau) = 0$ for $\tau < 0$. 

\begin{rem}
Even if we restrict the presentation to parameters $\theta_1,\theta_2 \in[0,\pi/2]$, which correspond to the most meaningful configurations from a practical perspective, developments reported in this paper readily extend to the case $\theta_1,\theta_2 \in[0,\pi)$ provided $q$ in (\ref{eq: writting of tilde_q}) is selected sufficiently large positive so that (\ref{eq: inner product Af and f}) still holds, implying in particular $\lambda_n \geq 0$ for all $n \geq 1$. In this case, one merely needs to modify the change of variable formula (\ref{eq: change of variable}) by the following one: $w(t,x) = z(t,x) - \frac{x^\alpha}{c_{\theta_2} + \alpha s_{\theta_2}} u(t)$ where $\alpha > 1$ is fixed so that $c_{\theta_2} + \alpha s_{\theta_2} \neq 0$.
\end{rem}

The system output $y(t) \in\R$ is selected as the either left Dirichlet trace (in this case $\theta_1 \in (0,\pi/2]$)
\begin{equation}\label{eq: Dirichlet output}
y_D (t) = z(t,0)
\end{equation}
or left Neumann trace (in this case $\theta_1 \in [0,\pi/2)$)
\begin{equation}\label{eq: Neumann output}
y_N (t) = z_x(t,0) .
\end{equation}
Without loss of generality, let $q \in\mathcal{C}^0([0,1])$ and $q_c \in\R$ be such that
\begin{equation}\label{eq: writting of tilde_q}
\tilde{q}(x) = q(x) - q_c , \quad q(x) > 0  .
\end{equation}
This allows to consider the Sturm-Liouville operator $\mathcal{A}$ and its related properties as described in Subsection~\ref{subsec: prop sturm liouville operators}. In particular, since $q \geq 0$, the eingenvalues $\lambda_n$ of $\mathcal{A}$ are such that $\lambda_n \geq 0$. Note however that the actual modes of the reaction-diffusion PDE (\ref{eq: PDE}) are given by $-\lambda_n + q_c$, hence a finite number of them may be unstable.

\subsection{Spectral reduction}

In order to obtain an equivalent homogeneous representation of  (\ref{eq: PDE}), we define the change of variable
\begin{equation}\label{eq: change of variable}
w(t,x) = z(t,x) - \frac{x^2}{c_{\theta_2} + 2 s_{\theta_2}} u(t-h) .
\end{equation}
Introducing $v = \dot{u}$, we infer that
\begin{subequations}\label{eq: PDE Dirichlet - homogeneous}
\begin{align}
& \dot{u}(t) = v(t) \\
& w_t(t,x) = (p(x) w_x(t,x))_x - \tilde{q}(x) w(t,x) \\
& \phantom{w_t(t,x) =} \;  + a(x) u(t-h) + b(x) v(t-h) \\ 
& c_{\theta_1} w(t,0) - s_{\theta_1} w_x(t,0) = 0 \\
& c_{\theta_2} w(t,1) + s_{\theta_2} w_x(t,1) = 0 \\
& w(0,x) = w_0(x)
\end{align}
\end{subequations}
where $a(x) = \frac{1}{c_{\theta_2} + 2 s_{\theta_2}} \{ 2p(x) + 2xp'(x) - x^2 \tilde{q}(x) \}$, $b(x) = -\frac{x^2}{c_{\theta_2} + 2 s_{\theta_2}}$, and $w_0(x) = z_0(x) - \frac{x^2}{c_{\theta_2} + 2 s_{\theta_2}} u(-h) = z_0(x)$. Let the coefficients of projection be defined by $z_n(t) = \left< z(t,\cdot) , \phi_n \right>$, $w_n(t) = \left< w(t,\cdot) , \phi_n \right>$, $a_n = \left< a , \phi_n \right>$, and $b_n = \left< b , \phi_n \right>$. In particular we have from (\ref{eq: change of variable}) that
\begin{equation}\label{eq: link z_n and w_n}
w_n(t) = z_n(t) + b_n u(t-h), \quad n \geq 1 .
\end{equation}
Using standard arguments, see e.g.~\cite{lhachemi2020boundary,mironchenko2020local} for details, the projection of (\ref{eq: PDE Dirichlet - homogeneous}) into the Hilbert basis $(\phi_n)_{n \geq 1}$ gives
\begin{subequations}\label{eq: dynamics w_n}
\begin{align}
\dot{u}(t) & = v(t) \\
\dot{w}_n(t) & = (-\lambda_n + q_c) w_n(t) + a_n u(t-h) + b_n v(t-h) 
\end{align}
\end{subequations}
with $w(t,\cdot) = \sum_{n \geq 1} w_n(t) \phi_n$ in $L^2$ norm for mild solutions and in $H^2$ norm for classical solutions (see the end of Section~\ref{subsec: prop sturm liouville operators}). Using (\ref{eq: link z_n and w_n}) into the latter idendity, the projection of (\ref{eq: PDE}) reads
\begin{equation}\label{eq: dynamics z_n}
\dot{z}_n(t) = (-\lambda_n + q_c) z_n(t) + \beta_n u(t-h)
\end{equation}
with $\beta_n = a_n + (-\lambda_n+q_c)b_n = p(1) \{ - c_{\theta_2} \phi_n'(1) + s_{\theta_2} \phi_n(1) \} = O(\sqrt{\lambda_n})$. Here we have $z(t,\cdot) = \sum_{n \geq 1} z_n(t) \phi_n$ in $L^2$ norm. Finally, when dealing with classical solutions, the Dirichlet measurement $y_D(t)$ given by (\ref{eq: Dirichlet output}) can be expressed as the series expansion:
\begin{equation}\label{eq: dynamics w_n - y Dirichlet}
y_D(t) = z(t,0) = w(t,0) = \sum_{n \geq 1} w_n(t) \phi_n(0) 
\end{equation}
while for the Neumann measurement $y_N(t)$ given by (\ref{eq: Neumann output}) we have:
\begin{equation}\label{eq: dynamics w_n - y Neumann}
y_N(t) = z_x(t,0) = w_x(t,0) = \sum_{n \geq 1} w_n(t) \phi_n'(0)  .
\end{equation}

\begin{rem}\label{rem: series expansion measurement}
Note that the series expansions (\ref{eq: dynamics w_n - y Dirichlet}-\ref{eq: dynamics w_n - y Neumann}) hold for the coefficients of projection $w_n$, i.e., for the PDE in $w$ coordinates. Such a series expansion does not hold for the coefficients of projection $z_n$, i.e., for the PDE in $z$ coordinates.
\end{rem}

\begin{rem}\label{rem: position of this work}
The use of a predictor feedback to achieve the boundary stabilization of (\ref{eq: PDE}) in the case of Neumann actuation and boundary condition ($\theta_1 = \theta_2 = \pi/2$) was reported first in~\cite{katz2021sub} for a bounded output operator, namely $y(t) = \int_0^1 c(x) z(t,x) \,\mathrm{d}x$ with $c \in L^2(0,1)$, and for system trajectories evaluated in $L^2$ norm. In this very specific setting, the authors managed to perform the both control design and stability analysis on the sole representation (\ref{eq: dynamics z_n}), i.e., for the PDE in original coordinates (\ref{eq: PDE}). This approach is strongly tailored for the above-mentioned setting and is hardly extendable to other types of boundary control (Dirichlet/Robin) and to unbounded measurement operators (Dirichlet/Neumann). This is essentially because the consideration of a Neumann actuation yields the most favorable case $\beta_n = p(1) \phi_n(1) = O(1)$ while any other boundary actuation setting (Dirichlet/Robin) gives in general no better than $\beta_n = O(\sqrt{\lambda_n})$. However, one of the crucial points of the $L^2$ stability analysis performed in~\cite{katz2021sub} relies in the use of the estimate $2 \sum \beta_n z_n u \leq \alpha \sum (\beta_n^2 / \lambda_n) u^2 + \alpha^{-1} \sum \lambda_n z_n^2$, valid for any $\alpha > 0$, where the convergence of the first series on the RHS holds for the Neumann actuation setting ($\beta_n = O(1)$) and with a term $\sum \lambda_n z_n^2$ that can be handled, provided its convergence, in the $L^2$ Lyapunov stability analysis\footnote{Essentially because $\lambda_n$ appearing in $\sum \lambda_n z_n^2 = \sum \lambda_n^\alpha z_n^2$ has a power $\alpha = 1$ that is not larger than the one in the dynamics of the modes (\ref{eq: dynamics z_n}).}. This approach fails in the case of Dirichlet/Robin actuation settings ($\beta_n = O(\sqrt{\lambda_n})$). The situation is even more stringent when trying to assess the stability of the system trajectories in $H^1$ norm (possibly with unbounded output operators instead of a bounded one) since this would led to a term of the form $\sum \lambda_n \beta_n z_n$ that cannot be neither handled with the above approach. In this paper, we completely solve the control design problem for the general reaction diffusion (\ref{eq: PDE}) with Dirichlet/Neumann/Robin boundary control/condition and with a measurement selected either as the Dirichlet (\ref{eq: Dirichlet output}) or Neumann (\ref{eq: Neumann output}) trace. The proposed strategy consists in designing the predictor feedback based on the representation (\ref{eq: dynamics z_n}) in original coordinates (\ref{eq: PDE}) while the stability analysis is performed based on (\ref{eq: dynamics w_n}) in homogeneous coordinates (\ref{eq: PDE Dirichlet - homogeneous}).
\end{rem}

\section{Case of a Dirichlet measurement}\label{sec: Case of a Dirichlet measurement}

We consider in this section the input-delayed reaction-diffusion system (\ref{eq: PDE}) for $\theta_1 \in (0,\pi/2]$ with Dirichlet measurement (\ref{eq: Dirichlet output}).

\subsection{Control strategy}\label{subsec: control strategy}

Let $\delta > 0$ be the desired exponential decay rate for the closed-loop system trajectories. Let $N_0 \geq 1$ be such that $-\lambda_n + q_c < - \delta < 0$ for all $n \geq N_0 + 1$. Let $N \geq N_0 + 1$ be arbitrarily given and that will be specified later. Consider first the following observer dynamics used to estimate the $N$ first modes of the plant in $z$-coordinates: 
\begin{subequations}\label{eq: controller part 1 - Dirichlet}
\begin{align}
\hat{w}_n(t) & = \hat{z}_n(t) + b_n u(t-h) \label{eq: controller 1 - Dirichlet} \\
\dot{\hat{z}}_n(t) & = (-\lambda_n+q_c) \hat{z}_n(t) + \beta_n u(t-h) \label{eq: controller 2 - Dirichlet} \\
& \phantom{=}\; - l_n \left\{ \sum_{k = 1}^N \hat{w}_k(t) \phi_k(0) - y_D(t) \right\}  ,\; 1 \leq n \leq N_0 \nonumber \\
\dot{\hat{z}}_n(t) & = (-\lambda_n+q_c) \hat{z}_n(t) + \beta_n u(t-h) ,\; N_0+1 \leq n \leq N \label{eq: controller 3 - Dirichlet}
\end{align}
\end{subequations}
where $l_n \in\R$ are the observer gains. 

\begin{rem}
Dynamics (\ref{eq: controller part 1 - Dirichlet}) constitutes an observer of the $N$ first modes $z_n$ of the PDE in (original) $z$ coordinates. However, due to Remark~\ref{rem: series expansion measurement} and in view of (\ref{eq: dynamics w_n - y Dirichlet}), the estimation $\hat{y}_D(t)$ of the actual Dirichlet measurement $y_D(t)$ is expressed in function of the estimation $\hat{w}_n$ of the modes $w_n$ of the PDE in homogeneous coordinates $w$ as $\hat{y}_D(t) = \sum_{k = 1}^N \hat{w}_k(t) \phi_k(0)$. Hence, even if (\ref{eq: controller part 1 - Dirichlet}) estimates the modes $z_n$ in $z$ coordinates, the correction of the error of measurement is done based on the modes $w_n$ in $w$ coordinates.
\end{rem}

\begin{rem}
The idea to split the observer dynamics into two parts, one with active correction of the estimation error for the first modes as in (\ref{eq: controller 2 - Dirichlet}) and one without correction of the estimation error as in (\ref{eq: controller 3 - Dirichlet}), roots back to \cite{sakawa1983feedback} in a delay-free context with bounded input and output operators. Since then, such an idea, sometimes referred to as the add of a ``Residual Mode Filter'' and which was shown to be of paramount importance for ensuring closed-loop stability~\cite{balas1988finite} has been extended in various directions~\cite{balas1988finite,grune2021finite,harkort2011finite,katz2020constructive,lhachemi2020finitePI,lhachemi2020finite,lhachemi2021local,lhachemi2021nonlinear,sakawa1983feedback,sano2012stability}. 
\end{rem}

Due to the input delay $h > 0$, we need to introduce a predictor component. To do so, let $\hat{Z}^{N_0} = \begin{bmatrix} \hat{z}_1 & \ldots & \hat{z}_{N_0} \end{bmatrix}^\top$, $A_0 = \mathrm{diag}(-\lambda_1 + q_c , \ldots , -\lambda_{N_0} + q_c)$, and $\mathfrak{B}_0 = \begin{bmatrix} \beta_1 & \ldots & \beta_{N_0} \end{bmatrix}^\top$. We can now introduce the following Artstein transformation~\cite{artstein1982linear}:
\begin{equation}\label{eq: Artstein transformation}
\hat{Z}_A^{N_0}(t) = e^{A_0 h} \hat{Z}^{N_0}(t) + \int_{t-h}^t e^{A_0 (t-\tau)} \mathfrak{B}_0 u(\tau) \,\mathrm{d}\tau .
\end{equation}
Then the control is defined as the predictor feedback~\cite{karafyllis2017predictor}:
\begin{equation}\label{eq: input Dirichlet}
u(t) = K \hat{Z}_A^{N_0}(t) , \quad t \geq 0 
\end{equation} 
where $K \in\R^{1 \times N_0}$ is the feedback gain.

\begin{rem}\label{rem WP1}
The well-posedness of the closed-loop system composed of the plant (\ref{eq: PDE Dirichlet - homogeneous}) in homogeneous coordinates $w$, the Dirichlet measurement (\ref{eq: Dirichlet output}), and the controller (\ref{eq: controller part 1 - Dirichlet}-\ref{eq: input Dirichlet}) in terms of classical solutions for initial condition $z_0 \in H^2(0,1)$ so that $c_{\theta_1} z_0(0) - s_{\theta_1} z_0'(0) = 0$ and $c_{\theta_2} z_0(1) + s_{\theta_2} z_0'(1) = 0$,  with null control in negative times ($u(\tau)=0$ for $\tau < 0$) and zero initial condition for the observer ($\hat{z}_n(0)=0$), is a direct consequence of~\cite[Thm.~6.3.1 and~6.3.3]{pazy2012semigroups} and the invertibility of the Artstein transformation~\cite{bresch2018new} using a classical induction argument. Having obtained classical solutions based on the homogeneous representation (\ref{eq: PDE Dirichlet - homogeneous}), standard arguments~\cite[Sec.~3.3]{curtain2012introduction} using the change of variable formula (\ref{eq: change of variable}) give the existence of classical solutions for the closed-loop system with the plant (\ref{eq: PDE}) expressed in original $z$ coordinates.
\end{rem}

In preparation of the statement of the main result, we define the matrices $A_1 = \mathrm{diag}(-\lambda_{N_0+1} + q_c , \ldots , -\lambda_{N} + q_c)$, $\tilde{\mathfrak{B}}_1 = \begin{bmatrix} \beta_{N_0 +1}/\lambda_{N_0 +1} & \ldots & \beta_N/\lambda_N \end{bmatrix}^\top$, $C_0 = \begin{bmatrix} \phi_1(0) & \ldots & \phi_{N_0}(0) \end{bmatrix}$, $\tilde{C}_1 = \begin{bmatrix} \frac{\phi_{N_0 +1}(0)}{\sqrt{\lambda_{N_0 +1}}} & \ldots & \frac{\phi_{N}(0)}{\sqrt{\lambda_{N}}} \end{bmatrix}$, $L = \begin{bmatrix} l_1 & \ldots & l_{N_0} \end{bmatrix}^\top$, 
\begin{equation*}
F =
\begin{bmatrix}
A_0 + \mathfrak{B}_0 K & e^{A_0 h} LC_0 & 0 & e^{A_0 h} L\tilde{C}_1 \\
0 & A_0 - L C_0 & 0 & - L\tilde{C}_1 \\
\tilde{\mathfrak{B}}_1 K & 0 & A_1 & 0 \\
0 & 0 & 0 & A_1
\end{bmatrix} ,\;
\mathcal{L} =
\begin{bmatrix}
e^{A_0 h} L \\ -L \\ 0 \\ 0
\end{bmatrix}
\end{equation*}
$E = \begin{bmatrix} A_0 + \mathfrak{B}_0 K & e^{A_0 h} LC_0 & 0 & e^{A_0 h} L\tilde{C}_1 & e^{A_0 h} L \end{bmatrix}$, and $\tilde{K} = \begin{bmatrix} K & 0 & 0 & 0 \end{bmatrix}$.

\begin{rem}\label{rem: kalman condition}
Since $A_0$ is diagonal, the Hautus test combined with the boundary conditions involved in the definition of $D(\mathcal{A})$ show that the pairs $(A_0,\mathfrak{B}_0)$ and $(A_0,C_0)$ both satisfy the Kalman condition. The same remark applies in the case of the Neumann boundary measurement (\ref{eq: Neumann output}) studied in Section~\ref{sec: Case of a Neumann measurement}. Subsequently, the feedback gain $K$ from (\ref{eq: input Dirichlet}) and the observer gain $L$ whose coefficients $l_n$ appear in (\ref{eq: controller part 1 - Dirichlet}) are computed so that  $A_0 + \mathfrak{B}_0 K$ and $A_0 - L C_0$ are Hurwitz with eigenvalues that have a real part strictly less than $-\delta<0$. To complete the control design procedure, it merely remains to select adequately the dimension $N$ of the observer to ensure the stability of the closed-loop system with exponential decay rate $\delta > 0$.
\end{rem}

\subsection{Main stability results}

\begin{thm}\label{thm1}
Let $\theta_1 \in (0,\pi/2]$, $\theta_2 \in [0,\pi/2]$, $p \in\mathcal{C}^2([0,1])$ with $p > 0$, and $\tilde{q} \in\mathcal{C}^0([0,1])$. Let $q \in\mathcal{C}^0([0,1])$ and $q_c \in\R$ be such that (\ref{eq: writting of tilde_q}) holds. Let $\delta > 0$ and $N_0 \geq 1$ be such that $-\lambda_n + q_c < - \delta$ for all $n \geq N_0 + 1$. Let $K\in\R^{1 \times N_0}$ and $L\in\R^{N_0}$ be such that $A_0 + \mathfrak{B}_0 K$ and $A_0 - L C_0$ are Hurwitz with eigenvalues that have a real part strictly less than $-\delta<0$. Let $h > 0$ be given. For a given $N \geq N_0 +1$, assume that there exist $P \succ 0$, $Q_1,Q_2 \succeq 0$, $\alpha>1$, and $\beta,\gamma > 0$ such that 
\begin{equation}\label{eq: thm1 - constraints}
\Theta_1 \preceq 0 ,\quad \Theta_2 \leq 0, \quad R_1 \preceq 0, \quad R_2 \preceq 0
\end{equation}
where
\begin{subequations}
\begin{align}
\Theta_1 & = \begin{bmatrix} F^\top P + P F + 2 \delta P + \tilde{Q}_1 & P \mathcal{L} \\ \mathcal{L}^\top P & -\beta \end{bmatrix} + E^\top Q_2 E \label{eq: def theta1 Dirichlet} \\
\Theta_2 & = 2\gamma\left\{ - \left( 1 - \frac{1}{\alpha} \right) \lambda_{N+1}+q_c + \delta \right\} + \beta M_\phi \\
R_1 & = -e^{-2\delta h} Q_1 + \alpha\gamma \Vert \mathcal{R}_N a \Vert_{L^2}^2 K^\top K \label{eq: def R1 Dirichlet} \\
R_2 & = -e^{-2\delta h} Q_2 + \alpha\gamma \Vert \mathcal{R}_N b \Vert_{L^2}^2 K^\top K \label{eq: def R2 Dirichlet}
\end{align}
\end{subequations}
with $\tilde{Q}_1 = \mathrm{diag}(Q_1,0,0,0)$ and $M_\phi = \sum_{n \geq N+1} \frac{\vert \phi_n(0) \vert^2}{\lambda_n} < +\infty$. Then there exists a constant $M > 0$ such that for any initial condition $z_0 \in H^2(0,1)$ so that $c_{\theta_1} z_0(0) - s_{\theta_1} z_0'(0) = 0$ and $c_{\theta_2} z_0(1) + s_{\theta_2} z_0'(1) = 0$, the trajectories of the closed-loop system composed of the plant (\ref{eq: PDE}), the Dirichlet measurement (\ref{eq: Dirichlet output}), and the controller (\ref{eq: controller part 1 - Dirichlet}-\ref{eq: input Dirichlet}) with null control in negative times ($u(\tau)=0$ for $\tau < 0$) and zero initial condition for the observer ($\hat{z}_n(0)=0$) satisfy $\Vert z(t,\cdot) \Vert_{H^1}^2 + \sup_{\tau\in[t-h,t]} \vert u(\tau) \vert^2 + \sum_{n = 1}^{N} \hat{z}_n(t)^2 \leq M e^{-2\delta t} \Vert z_0 \Vert_{H^1}^2$ for all $t \geq 0$. Moreover, for any given $h > 0$, the constraints (\ref{eq: thm1 - constraints}) are always feasible for $N$ selected large enough.
\end{thm}

\textbf{Proof.} 
We first write a finite dimensional model capturing the $N$ first modes of the PDE and the dynamics (\ref{eq: controller part 1 - Dirichlet}-\ref{eq: input Dirichlet}) of the output feedback controller. We define $e_n = z_n - \hat{z}_n$ for all $1 \leq n \leq N$. In view of (\ref{eq: controller 1 - Dirichlet}-\ref{eq: controller 2 - Dirichlet}) and based on (\ref{eq: link z_n and w_n}) and (\ref{eq: dynamics w_n - y Dirichlet}), we obtain that
\begin{equation}\label{eq: controller 2 bis}
\dot{\hat{z}}_n 
= (-\lambda_n + q_c) \hat{z}_n + \beta_n u(\cdot-h) + l_n \sum_{k=1}^{N} \phi_k(0) e_k + l_n \zeta
\end{equation}
for $1 \leq n \leq N_0$ where $\zeta(t) = \sum_{n \geq N+1} w_n(t) \phi_n(0)$. Introducing $E^{N_0} = \begin{bmatrix} e_1 & \ldots & e_{N_0} \end{bmatrix}^\top$, the scaled error $\tilde{e}_n = \sqrt{\lambda_n} e_n$ (see~\cite{lhachemi2020finite}), and $\tilde{E}^{N - N_0} = \begin{bmatrix} \tilde{e}_{N_0 +1} & \ldots & \tilde{e}_{N} \end{bmatrix}^\top$, we deduce that
\begin{equation}\label{eq: time derivative hat_Z_N0}
\dot{\hat{Z}}^{N_0} = A_0 \hat{Z}^{N_0} + \mathfrak{B}_0 u(\cdot-h) + LC_0 E^{N_0} + L\tilde{C}_1 \tilde{E}^{N-N_0} + L \zeta .
\end{equation}
Taking first the time derivative of (\ref{eq: Artstein transformation}) and then inserting (\ref{eq: input Dirichlet}) and (\ref{eq: time derivative hat_Z_N0}), we infer that
\begin{align}
\dot{\hat{Z}}_A^{N_0} & = ( A_0 + \mathfrak{B}_0 K ) \hat{Z}_A^{N_0} \label{eq: truncated model - 4 ODEs - 1} \\
& \phantom{=}\; + e^{A_0 h} \left( LC_0 E^{N_0} + L\tilde{C}_1 \tilde{E}^{N-N_0} + L \zeta \right) . \nonumber
\end{align}
Consider now (\ref{eq: controller 3 - Dirichlet}). Defining the scaled estimation $\tilde{z}_n = \hat{z}_n / \lambda_n$ and $\tilde{Z}^{N-N_0} = \begin{bmatrix} \tilde{z}_{N_0 + 1} & \ldots & \tilde{z}_{N} \end{bmatrix}^\top$, this latter dynamics can be written as
\begin{equation*}
\dot{\tilde{Z}}^{N-N_0} = A_1 \tilde{Z}^{N-N_0} + \tilde{\mathfrak{B}}_1 u(t-h) .
\end{equation*}
In order to eliminate the input delay, we consider the second following Artstein transformation:
\begin{equation}\label{eq: Artstein transformation bis}
\tilde{Z}_A^{N-N_0}(t) = e^{A_1 h} \tilde{Z}^{N-N_0}(t) + \int_{t-h}^t e^{A_1 (t-\tau)} \tilde{\mathfrak{B}}_1 u(\tau) \,\mathrm{d}\tau
\end{equation}
which implies, along with (\ref{eq: input Dirichlet}), that
\begin{equation}\label{eq: truncated model - 4 ODEs - 2}
\dot{\tilde{Z}}_A^{N-N_0} = A_1 \tilde{Z}_A^{N-N_0} + \tilde{\mathfrak{B}}_1 K \hat{Z}_A^{N_0} .
\end{equation}
Combining (\ref{eq: dynamics z_n}) and (\ref{eq: controller 2 - Dirichlet}-\ref{eq: controller 3 - Dirichlet}), the error dynamics reads
\begin{subequations}\label{eq: truncated model - 4 ODEs - 3}
\begin{align}
\dot{E}^{N_0} & = ( A_0 - L C_0 ) E^{N_0} - L \tilde{C}_1 \tilde{E}^{N-N_0} - L \zeta , \\
\dot{\tilde{E}}^{N-N_0} & = A_1 \tilde{E}^{N-N_0} .
\end{align}
\end{subequations}
Therefore, defining the vector
\begin{equation}\label{eq: truncated model - def X}
X = \mathrm{col}\left( \hat{Z}_A^{N_0} , E^{N_0} , \tilde{Z}_A^{N-N_0} , \tilde{E}^{N-N_0} \right) ,
\end{equation}
we infer from (\ref{eq: truncated model - 4 ODEs - 1}) and (\ref{eq: truncated model - 4 ODEs - 2}-\ref{eq: truncated model - 4 ODEs - 3}) that
\begin{equation}\label{eq: truncated model}
\dot{X} = F X + \mathcal{L} \zeta .
\end{equation}
Finally, defining $\tilde{X} = \mathrm{col}\left( X , \zeta \right)$ and based on (\ref{eq: input Dirichlet}) and (\ref{eq: truncated model - 4 ODEs - 1}), we also have 
\begin{equation}\label{eq: derivative v of command input u}
u = \tilde{K} X  , \quad  \quad v = \dot{u} = K \dot{\hat{Z}}_A^{N_0} , \quad \dot{\hat{Z}}_A^{N_0} = E \tilde{X} .
\end{equation}

We can now perform the stability analysis. Consider the functional defined by $V(t) = V_0(t)+V_1(t)+V_2(t)$ where
\begin{subequations}\label{eq: Lyapunov function H1 norm}
\begin{align}
V_0(t) & = X(t)^\top P X(t) + \gamma \sum_{n \geq N+1} \lambda_n w_n(t)^2 \\
V_1(t) & = \int_{(t-h)_+}^t e^{-2 \delta (t-s)} \hat{Z}_A^{N_0}(\tau)^\top Q_1 \hat{Z}_A^{N_0}(\tau) \mathrm{d}\tau \\
V_2(t) & = \int_{(t-h)_+}^t e^{-2 \delta (t-s)} \dot{\hat{Z}}_A^{N_0}(\tau)^\top Q_2 \dot{\hat{Z}}_A^{N_0}(\tau) \,\mathrm{d}\tau
\end{align}
\end{subequations}
with $(t-h)_+ = \max(t-h,0)$. The computation of the time derivative of $V$ for $t > h$ gives 
\begin{align*}
& \dot{V} \leq \tilde{X}^\top \begin{bmatrix} F^\top P + P F & P \mathcal{L} \\ \mathcal{L}^\top P & 0 \end{bmatrix} \tilde{X} + 2\gamma \sum_{n \geq N+1} \lambda_n (-\lambda_n + q_c) w_n^2 \\
& + 2\gamma \sum_{n \geq N+1} \lambda_n \{ a_n u(\cdot-h) + b_n v(\cdot-h) \} w_n - 2\delta V_1 - 2\delta V_2 \\
& + (\hat{Z}_A^{N_0})^\top Q_1 \hat{Z}_A^{N_0} - e^{-2\delta h} \hat{Z}_A^{N_0}(\cdot-h)^\top Q_1 \hat{Z}_A^{N_0}(\cdot-h) \\
& + (\dot{\hat{Z}}_A^{N_0})^\top Q_2 \dot{\hat{Z}}_A^{N_0} - e^{-2\delta h} \dot{\hat{Z}}_A^{N_0}(\cdot-h)^\top Q_2 \dot{\hat{Z}}_A^{N_0}(\cdot-h) .
\end{align*}
Using Young inequality and invoking (\ref{eq: input Dirichlet}), we obtain that
\begin{align*}
& 2 \sum_{n \geq N+1} \lambda_n a_n u(\cdot-h) w_n \\
& \leq \frac{1}{\alpha} \sum_{n \geq N+1} \lambda_n^2 w_n^2
+ \alpha \Vert \mathcal{R}_N a \Vert_{L^2}^2 u(\cdot-h)^2 \\
& \leq \frac{1}{\alpha} \sum_{n \geq N+1} \lambda_n^2 w_n^2
+ \alpha \Vert \mathcal{R}_N a \Vert_{L^2}^2 \hat{Z}_A^{N_0}(\cdot-h)^\top K^\top K \hat{Z}_A^{N_0}(\cdot-h)
\end{align*}
for any $\alpha > 0$ and, similarly,
\begin{align*}
& 2 \sum_{n \geq N+1} \lambda_n b_n v(\cdot-h) w_n \\
& \leq \frac{1}{\alpha} \sum_{n \geq N+1} \lambda_n^2 w_n^2
+ \alpha \Vert \mathcal{R}_N b \Vert_{L^2}^2 \dot{\hat{Z}}_A^{N_0}(\cdot-h)^\top K^\top K \dot{\hat{Z}}_A^{N_0}(\cdot-h)
\end{align*}
Combining the latter estimates and using (\ref{eq: derivative v of command input u}), we obtain
\begin{align*}
& \dot{V} + 2 \delta V \\
& \leq \tilde{X}^\top \left\{ \begin{bmatrix} F^\top P + P F + 2 \delta P + \tilde{Q}_1 & P \mathcal{L} \\ \mathcal{L}^\top P & 0 \end{bmatrix} + E^\top Q_2 E \right\} \tilde{X} \\
& + 2\gamma \sum_{n \geq N+1} \lambda_n \left\{ - \left( 1 - \frac{1}{\alpha} \right) \lambda_n + q_c + \delta \right\} w_n^2 \\
& + \hat{Z}_A^{N_0}(\cdot-h)^\top R_1 \hat{Z}_A^{N_0}(\cdot-h) + \dot{\hat{Z}}_A^{N_0}(\cdot-h)^\top R_2 \dot{\hat{Z}}_A^{N_0}(\cdot-h) .
\end{align*}
Recalling that $\zeta = \sum_{n \geq N+1} w_n \phi_n(0)$ we have $\zeta^2 \leq M_\phi \sum_{n \geq N+1} \lambda_n w_n^2$. Hence, we obtain, for any $\beta > 0$,
\begin{align}
& \dot{V} + 2 \delta V \leq \tilde{X}^\top \Theta_1 \tilde{X} + \sum_{n \geq N+1} \lambda_n \Gamma_n w_n^2 \nonumber \\
& + \hat{Z}_A^{N_0}(\cdot-h)^\top R_1 \hat{Z}_A^{N_0}(\cdot-h) + \dot{\hat{Z}}_A^{N_0}(\cdot-h)^\top R_2 \dot{\hat{Z}}_A^{N_0}(\cdot-h) \label{eq: estimate dotV+deltaV}
\end{align}
where $\Gamma_n = 2 \gamma \left\{ - \left( 1 - \frac{1}{\alpha} \right) \lambda_n + q_c + \delta \right\} + \beta M_\phi$. For $\alpha > 1$, we have $\Gamma_n \leq \Gamma_{N+1} = \Theta_2$ for all $n \geq N+1$. Thus we infer from (\ref{eq: thm1 - constraints}) that $\dot{V} + 2 \delta V \leq 0$ for all $t > h$, implying that $V(t) \leq e^{-2\delta(t-h)}V(h)$ for all $t \geq h$. 

Computing now the time derivative of $V$ for $t \in (0,h)$ for which $u(\cdot-h)$ and $v(\cdot-h)$ are zero, and proceeding similarly to the previous paragraph, we infer that  $\dot{V}+2\delta V \leq \tilde{X}^\top \Theta_1 \tilde{X} + \sum_{n \geq N+1} \lambda_n \Gamma_n' w_n^2$ with $\Gamma_n' = 2 \gamma \left\{ - \lambda_n + q_c + \delta \right\} + \beta M_\phi \leq \Gamma_n \leq \Theta_2$ for $n \geq N+1$. Hence we obtain that $\dot{V}+2\delta V \leq 0$ on $(0,h)$ thus $V(t) \leq e^{-2\delta t}V(0)$ for all $0 \leq t \leq h$. Combining this estimate with the result of the previous paragraph, we infer that $V(t) \leq e^{-2\delta t} V(0)$ for all $t \geq 0$. The claimed stability estimate now easily follows from the definition of $V$, the estimates (\ref{eq: inner product Af and f}), the control law (\ref{eq: input Dirichlet}), and the two Artstein transformations (\ref{eq: Artstein transformation}) and (\ref{eq: Artstein transformation bis}).

It remains to show that the constraints $\Theta_1,\Theta_2,R_1,R_2 \preceq 0$ are feasible when selecting $N \geq N_0 +1$ large enough. Regarding the matrix $F + \delta I$, we note that (i) $A_0 + \mathfrak{B}_0 K + \delta I$ and $A_0 - L C_0 + \delta I$ are Hurwitz; (ii) $\Vert e^{(A_1 + \delta I)t} \Vert \leq e^{- \kappa_0 t}$ for all $t \ge 0$ with $\kappa_0 = \lambda_{N_0+1} - q_c - \delta > 0$ is independent of $N$; and (iii) $\Vert e^{A_0 h} L\tilde{C}_1 \Vert \leq \Vert e^{A_0 h} \Vert \Vert L \Vert \Vert \tilde{C}_1 \Vert$, $\Vert L\tilde{C}_1 \Vert \leq \Vert L \Vert \Vert \tilde{C}_1 \Vert$ and $\Vert \tilde{\mathfrak{B}}_1 K \Vert \leq \Vert \tilde{\mathfrak{B}}_1 \Vert \Vert K \Vert$ where $e^{A_0 h}$, $K$, and $L$ are independent of the number of observed modes $N$ while $\Vert \tilde{C}_1 \Vert = O(1)$ and $\Vert \tilde{\mathfrak{B}}_1 \Vert = O(1)$ when $N \rightarrow +\infty$. Hence, applying the Lemma reported in Appendix to the matrix $F+\delta I$, we infer that the solution $P \succ 0$ to $F^\top P + P F + 2 \delta P = - I$ is such that $\Vert P \Vert = O(1)$ as $N \rightarrow + \infty$. Note also that $\Vert \mathcal{L} \Vert$ is a constant independent of $N$ while $M_\phi = O(1)$ and $\Vert E \Vert = O(1)$ as $N \rightarrow + \infty$. We fix arbitrarily the value of $\alpha > 1$ and we set $\beta = \sqrt{N}$, $\gamma = 1/N$, $Q_1 = 2 e^{2\delta h} \alpha\gamma \Vert \mathcal{R}_N a \Vert_{L^2}^2 K^\top K$, and $Q_2 = 2 e^{2\delta h} \alpha\gamma \Vert \mathcal{R}_N b \Vert_{L^2}^2 K^\top K$. Hence, using in particular Schur complement, we infer that (\ref{eq: thm1 - constraints}) hold for $N \geq N_0 + 1$ selected large enough.
\qed

\begin{rem}
Recall the definition of $V(t) = V_0(t)+V_1(t)+V_2(t)$ with $V_i$ defined by (\ref{eq: Lyapunov function H1 norm}). The term $V_0$ presents a first term that accounts for the dynamics of the observer as well as the dynamics of the $N$ first modes of the (original) PDE in $z$ coordinates. However, the second term, which accounts for the residual modes $n \geq N+1$, is expressed in homogeneous $w$ coordinates. This point is key in the success of the stability assessment of the previous theorem.  The term $V_1$ is introduced as in~\cite{katz2021sub} in order to compensate the term $u(\cdot-h)$ appearing in the time derivative of $w_n$, see (\ref{eq: dynamics w_n}). Moreover, since the modes $n \geq N+1$ are captured by $w_n$ in $w$ coordinates, the time derivative of $w_n$ also implies the occurrence of $v = \dot{u} = K \dot{\hat{Z}}_A^{N_0}$, see again (\ref{eq: dynamics w_n}). This latter term is handled in the stability analysis by the introduction of the term $V_2$.
\end{rem}

While Theorem~\ref{thm1} assesses the stability of the closed-loop in $H^1$ norm, the following theorem states a similar result but for trajectories evaluated in $L^2$ norm.

\begin{thm}\label{thm2}
Let $\theta_1 \in (0,\pi/2]$, $\theta_2 \in [0,\pi/2]$, $p \in\mathcal{C}^2([0,1])$ with $p > 0$, and $\tilde{q} \in\mathcal{C}^0([0,1])$. Let $q \in\mathcal{C}^0([0,1])$ and $q_c \in\R$ be such that (\ref{eq: writting of tilde_q}) holds. Let $\delta > 0$ and $N_0 \geq 1$ be such that $-\lambda_n + q_c < - \delta$ for all $n \geq N_0 + 1$. Let $K\in\R^{1 \times N_0}$ and $L\in\R^{N_0}$ be such that $A_0 + \mathfrak{B}_0 K$ and $A_0 - L C_0$ are Hurwitz with eigenvalues that have a real part strictly less than $-\delta<0$. Let $h > 0$ be given. For a given $N \geq N_0 +1$, assume that there exist $P \succ 0$, $Q_1,Q_2 \succeq 0$, and $\alpha,\beta,\gamma > 0$ such that 
\begin{equation}\label{eq: thm2 - constraints}
\Theta_1 \preceq 0 ,\quad \Theta_2 \leq 0, \quad \Theta_3 \geq 0 , \quad R_1 \preceq 0, \quad R_2 \preceq 0
\end{equation}
where $\Theta_1,R_1,R_2$ are defined by (\ref{eq: def theta1 Dirichlet}), (\ref{eq: def R1 Dirichlet}), and (\ref{eq: def R2 Dirichlet}), respectively, while
\begin{align*}
\Theta_2 & = 2\gamma\left\{ - \lambda_{N+1} + q_c + \delta + \frac{1}{\alpha} \right\} + \beta M_\phi \lambda_{N+1}^{3/4} \\
\Theta_3 & = 2\gamma - \frac{\beta M_\phi}{\lambda_{N+1}^{1/4}}
\end{align*}
with $\tilde{Q}_1 = \mathrm{diag}(Q_1,0,0,0)$ and $M_\phi = \sum_{n \geq N+1} \frac{\vert \phi_n(0) \vert^2}{\lambda_n^{3/4}} < +\infty$. Then there exists a constant $M > 0$ such that for any initial condition $z_0 \in H^2(0,1)$ so that $c_{\theta_1} z_0(0) - s_{\theta_1} z_0'(0) = 0$ and $c_{\theta_2} z_0(1) + s_{\theta_2} z_0'(1) = 0$, the trajectories of the closed-loop system composed of the plant (\ref{eq: PDE}), the Dirichlet measurement (\ref{eq: Dirichlet output}), and the controller (\ref{eq: controller part 1 - Dirichlet}-\ref{eq: input Dirichlet}) with null control in negative times ($u(\tau)=0$ for $\tau < 0$) and zero initial condition for the observer ($\hat{z}_n(0)=0$) satisfy $\Vert z(t,\cdot) \Vert_{L^2}^2 + \sup_{\tau\in[t-h,t]} \vert u(\tau) \vert^2 + \sum_{n = 1}^{N} \hat{z}_n(t)^2 \leq M e^{-2\delta t} \Vert z_0 \Vert_{L^2}^2$ for all $t \geq 0$. Moreover, for any given $h > 0$, the constraints (\ref{eq: thm2 - constraints}) are always feasible for $N$ selected large enough.
\end{thm}

\textbf{Proof.} 
Consider the functional defined by $V(t) = V_0(t)+V_1(t)+V_2(t)$ where $V_0(t) = X(t)^\top P X(t) + \gamma \sum_{n \geq N+1} w_n^2$ while $V_1,V_2$ are defined as in (\ref{eq: Lyapunov function H1 norm}). Proceeding as in the proof of Theorem~\ref{thm1} but replacing the estimate of $\zeta$ by the following: $\zeta^2 \leq M_\phi \sum_{n \geq N+1} \lambda_n^{3/4} w_n^2$, we infer that
\begin{align*}
& \dot{V} + 2 \delta V \leq \tilde{X}^\top \Theta_1 \tilde{X} + \sum_{n \geq N+1} \Gamma_n w_n^2 \nonumber \\
& + \hat{Z}_A^{N_0}(\cdot-h)^\top R_1 \hat{Z}_A^{N_0}(\cdot-h) + \dot{\hat{Z}}_A^{N_0}(\cdot-h)^\top R_2 \dot{\hat{Z}}_A^{N_0}(\cdot-h)
\end{align*}
holds for $t > h$ with $\Gamma_n = 2 \gamma \left\{ - \lambda_n + q_c + \delta + \frac{1}{\alpha} \right\} + \beta M_\phi \lambda_{n}^{3/4}$. For $n \geq N+1$ we note that $\lambda_n^{3/4} = \lambda_n/\lambda_n^{1/4} \leq \lambda_n/\lambda_{N+1}^{1/4}$ hence 
\begin{align*}
\Gamma_n & \leq - \Theta_3 \lambda_n + 2\gamma \left\{ q_c + \delta + \frac{1}{\alpha} \right\} \\
& \leq - \Theta_3 \lambda_{N+1} + 2\gamma \left\{ q_c + \delta + \frac{1}{\alpha} \right\} = \Theta_2 \leq 0
\end{align*}
where we used that $\Theta_3 \geq 0$. Therefore, the assumptions imply that $\dot{V} + 2 \delta V \leq 0$ for $t > h$. Similarly to the proof of Theorem~\ref{thm1}, it can also be seen that $\dot{V} + 2 \delta V \leq 0$ for $t \in (0,h)$. Gathering these two results together, we infer that $V(t) \leq e^{-2\delta t} V(0)$ for all $t \geq 0$, implying the claimed stability estimate.

Regarding the feasibility of the constraints (\ref{eq: thm2 - constraints}) for $N \geq N_0 +1$ large enough, this can be achieved following the same procedure as in the proof of Theorem~\ref{thm1} with $\alpha > 0$ arbitrarily fixed, $\beta = N^{1/8}$, and $\gamma = 1/N^{1/4}$.
\qed

\begin{rem}\label{rem LMIs}
Let  $N \geq N_0 + 1$ be a given number of modes to be observed. When fixing the value of $\alpha > 1$ (resp. $\alpha > 0$), the constraints (\ref{eq: thm1 - constraints}) from Theorem~\ref{thm1} (resp. the constraints constraints (\ref{eq: thm2 - constraints}) from Theorem~\ref{thm2}) take the form of LMIs for which efficient solvers exist. Moreover, as shown in the proof of the two theorems, the resulting LMI constraints remain feasible (when fixing arbitrarily the value of $\alpha$) for $N$ selected large enough.
\end{rem}

\section{Case of a Neumann measurement}\label{sec: Case of a Neumann measurement}

We consider in this section the input-delayed reaction-diffusion system (\ref{eq: PDE}) for $\theta_1 \in [0,\pi/2)$ with Neumann measurement (\ref{eq: Neumann output}).

\subsection{Control strategy}

Let $\delta > 0$ and $N_0 \geq 1$ be such that $-\lambda_n + q_c < - \delta < 0$ for all $n \geq N_0 + 1$. Let $N \geq N_0 + 1$ be arbitrarily given. Consider first the following observer dynamics used to estimate the $N$ first modes of the plant in $z$-coordinates: 
\begin{subequations}\label{eq: controller part 1 - Neumann}
\begin{align}
\hat{w}_n(t) & = \hat{z}_n(t) + b_n u(t-h) \label{eq: controller 1 - Neumann} \\
\dot{\hat{z}}_n(t) & = (-\lambda_n+q_c) \hat{z}_n(t) + \beta_n u(t-h) \label{eq: controller 2 - Neumann} \\
& \phantom{=}\; - l_n \left\{ \sum_{k = 1}^N \hat{w}_k(t) \phi_k'(0) - y_N(t) \right\}  ,\; 1 \leq n \leq N_0 \nonumber \\
\dot{\hat{z}}_n(t) & = (-\lambda_n+q_c) \hat{z}_n(t) + \beta_n u(t-h) ,\; N_0+1 \leq n \leq N \label{eq: controller 3 - Neumann}
\end{align}
\end{subequations}
where $l_n \in\R$ are the observer gains. With the same notations that the ones of the previous section and introducing the Artstein transformation (\ref{eq: Artstein transformation}), the control input is defined as
\begin{equation}\label{eq: input Neumann}
u(t) = K \hat{Z}_A^{N_0}(t) , \quad t \geq 0 
\end{equation} 
where $K \in\R^{1 \times N_0}$ is the feedback gain.

\begin{rem}
The well-posedness of the closed-loop system in terms of classical solutions for initial condition $z_0 \in H^2(0,1)$ so that $c_{\theta_1} z_0(0) - s_{\theta_1} z_0'(0) = 0$ and $c_{\theta_2} z_0(1) + s_{\theta_2} z_0'(1) = 0$,  with null control in negative times ($u(\tau)=0$ for $\tau < 0$) and zero initial condition for the observer ($\hat{z}_n(0)=0$), follows the same arguments as Remark~\ref{rem WP1}.
\end{rem}

We finally introduce the same matrices as at the end of Subsection~\ref{subsec: control strategy} except that we replace the definitions of $C_0,\tilde{C}_1$ by $C_0 = \begin{bmatrix} \phi_1'(0) & \ldots & \phi_{N_0}'(0) \end{bmatrix}$ and $\tilde{C}_1 = \begin{bmatrix} \phi_{N_0 +1}'(0)/\lambda_{N_0 +1} & \ldots & \phi_{N}'(0)/\lambda_{N} \end{bmatrix}$.

\subsection{Main stability result}

\begin{thm}\label{thm3}
Let $\theta_1 \in [0,\pi/2)$, $\theta_2 \in [0,\pi/2]$, $p \in\mathcal{C}^2([0,1])$ with $p > 0$, and $\tilde{q} \in\mathcal{C}^0([0,1])$. Let $q \in\mathcal{C}^0([0,1])$ and $q_c \in\R$ be such that (\ref{eq: writting of tilde_q}) holds. Let $\delta > 0$ and $N_0 \geq 1$ be such that $-\lambda_n + q_c < - \delta$ for all $n \geq N_0 + 1$. Let $K\in\R^{1 \times N_0}$ and $L\in\R^{N_0}$ be such that $A_0 + \mathfrak{B}_0 K$ and $A_0 - L C_0$ are Hurwitz with eigenvalues that have a real part strictly less than $-\delta<0$. Let $h > 0$ be given. For a given $N \geq N_0 +1$, assume that there exist $P \succ 0$, $Q_1,Q_2 \succeq 0$, $\epsilon\in(0,1/2]$, $\alpha>1$, and $\beta,\gamma > 0$ such that 
\begin{equation}\label{eq: thm3 - constraints}
\Theta_1 \preceq 0 ,\quad \Theta_2 \leq 0, \quad \Theta_3 \geq 0 , \quad R_1 \preceq 0, \quad R_2 \preceq 0
\end{equation}
where $\Theta_1,R_1,R_2$ are defined by (\ref{eq: def theta1 Dirichlet}), (\ref{eq: def R1 Dirichlet}), and (\ref{eq: def R2 Dirichlet}), respectively, while
\begin{align*}
\Theta_2 & = 2\gamma\left\{ - \left( 1 - \frac{1}{\alpha} \right) \lambda_{N+1} + q_c + \delta \right\} + \beta M_\phi(\epsilon) \lambda_{N+1}^{1/2+\epsilon} \\
\Theta_3 & = 2\gamma\left(1-\dfrac{1}{\alpha}\right) - \frac{\beta M_\phi(\epsilon)}{\lambda_{N+1}^{1/2-\epsilon}}
\end{align*}
with $\tilde{Q}_1 = \mathrm{diag}(Q_1,0,0,0)$ and $M_\phi(\epsilon) = \sum_{n \geq N+1} \frac{\vert \phi_n'(0) \vert^2}{\lambda_n^{3/2+\epsilon}} < +\infty$. Then there exists a constant $M > 0$ such that for any initial condition $z_0 \in H^2(0,1)$ so that $c_{\theta_1} z_0(0) - s_{\theta_1} z_0'(0) = 0$ and $c_{\theta_2} z_0(1) + s_{\theta_2} z_0'(1) = 0$, the trajectories of the closed-loop system composed of the plant (\ref{eq: PDE}), the Neumann measurement (\ref{eq: Neumann output}), and the controller (\ref{eq: controller part 1 - Neumann}-\ref{eq: input Neumann}) with null control in negative times ($u(\tau)=0$ for $\tau < 0$) and zero initial condition for the observer ($\hat{z}_n(0)=0$) satisfy $\Vert z(t,\cdot) \Vert_{H^1}^2 + \sup_{\tau\in[t-h,t]} \vert u(\tau) \vert^2 + \sum_{n = 1}^{N} \hat{z}_n(t)^2 \leq M e^{-2\delta t} \Vert z_0 \Vert_{H^1}^2$ for all $t \geq 0$. Moreover, for any given $h > 0$, the constraints (\ref{eq: thm3 - constraints}) are always feasible for $N$ selected large enough.
\end{thm}

\textbf{Proof.}
Proceeding as in the first part of the proof of Theorem~\ref{thm1} but with $\tilde{e}_n = \lambda_n e_n$ and $\zeta(t) = \sum_{n \geq N+1} w_n(t) \phi_n'(0)$, we infer that (\ref{eq: truncated model}) holds.

Consider the functional defined by $V(t) = V_0(t)+V_1(t)+V_2(t)$ where $V_0,V_1,V_2$ are defined by (\ref{eq: Lyapunov function H1 norm}). Proceeding as in the proof of Theorem~\ref{thm1} but replacing the estimate of $\zeta$ by: $\zeta^2 \leq M_\phi(\epsilon) \sum_{n \geq N+1} \lambda_n^{3/2+\epsilon} w_n^2$, we infer that
\begin{align*}
& \dot{V} + 2 \delta V \leq \tilde{X}^\top \Theta_1 \tilde{X} + \sum_{n \geq N+1} \lambda_n\Gamma_n w_n^2 \nonumber \\
& + \hat{Z}_A^{N_0}(\cdot-h)^\top R_1 \hat{Z}_A^{N_0}(\cdot-h) + \dot{\hat{Z}}_A^{N_0}(\cdot-h)^\top R_2 \dot{\hat{Z}}_A^{N_0}(\cdot-h)
\end{align*}
holds for $t > h$ with $\Gamma_n = 2 \gamma \left\{ - (1-\frac{1}{\alpha}) \lambda_n + q_c + \delta \right\} + \beta M_\phi(\epsilon) \lambda_{n}^{1/2 + \epsilon}$. For $n \geq N+1$ we note that $\lambda_n^{1/2 + \epsilon} = \lambda_n/\lambda_n^{1/2 - \epsilon} \leq \lambda_n/\lambda_{N+1}^{1/2 - \epsilon}$ hence 
\begin{align*}
\Gamma_n & \leq - \Theta_3 \lambda_n + 2\gamma \left\{ q_c + \delta \right\} \\
& \leq - \Theta_3 \lambda_{N+1} + 2\gamma \left\{ q_c + \delta \right\} = \Theta_2 \leq 0
\end{align*}
where we used that $\Theta_3 \geq 0$. Therefore, the assumptions imply that $\dot{V} + 2 \delta V \leq 0$ for $t > h$. Similarly to the proof of Theorem~\ref{thm1}, it can also be seen that $\dot{V} + 2 \delta V \leq 0$ for $t \in (0,h)$. Gathering these two results together, we infer that $V(t) \leq e^{-2\delta t} V(0)$ for all $t \geq 0$, implying the claimed stability estimate.

Regarding the feasibility of the constraints (\ref{eq: thm3 - constraints}) for $N \geq N_0 +1$ large enough, this is achieved following the same procedure as in the proof of Theorem~\ref{thm1} with $\alpha > 1$ arbitrarily fixed, $\epsilon = 1/8$, $\beta = N^{1/8}$, and $\gamma = 1/N^{3/16}$.
\qed

\begin{rem}
Similarly to Remark~\ref{rem LMIs}, LMIs can be derived from the constraints (\ref{thm3}) of Theorem~\ref{thm3} by fixing the values of $\alpha > 1$ and $\epsilon\in(0,1/2]$. Moreover, as shown in the proof of the theorem, the subsequent LMI conditions with any fixed $\alpha > 1$ and when setting $\epsilon = 1/8$ remain feasible for a number of observed modes $N$ selected large enough.
\end{rem}

\begin{rem}
In the case of a Dirichlet measurement, it was possible to propose a $L^2$ version of the stability result, namely Theorem~\ref{thm2}. However, the approach used in the proof of this latter result fails when trying to study the trajectories in $L^2$ norm for a Neumann measurement. This is because, in this setting, $\zeta = \sum_{n \geq N+1} w_n \phi_n'(0)$ hence $\zeta^2 \leq M_\phi(\epsilon) \sum_{n \geq N+1} \lambda_n^{3/2+\epsilon} w_n^2$ with $M_\phi(\epsilon) = \sum_{n \geq N+1} \frac{\vert \phi_n'(0) \vert^2}{\lambda_n^{3/2+\epsilon}}$ where all the terms are finite provided $\epsilon \in (0,1/2]$. However a term in $\lambda_n^{3/2+\epsilon}$ cannot be asymptotically dominated by a term in $\lambda_n$. Hence the procedure of Theorem~\ref{thm2} for a Dirichlet measurement cannot be used anymore in the case of a Neumann measurement.  
\end{rem}

\section{Numerical illustration}\label{sec: numerical example}

We consider the parameters $p=1$, $\tilde{q}=-5$, $\theta_1 = \pi/5$, $\theta_2 = 0$ (Dirichlet boundary control), and the input delay $h = 1\,\mathrm{s}$. In this setting, the reaction-diffusion PDE described by (\ref{eq: PDE}) is open-loop unstable. We set the feedback gain $K = -0.6950$. The observer gain is set as $L = 1.7695$ in the case of the Dirichlet measurement (\ref{eq: Dirichlet output}) while $L = 1.2856$ in the case of the Neumann measurement (\ref{eq: Neumann output}). With $\delta = 0.5$ and for the Dirichlet measurement (\ref{eq: Dirichlet output}), the constraints of Theorems~\ref{thm1} and~\ref{thm2} are feasible for $N = 2$ modes estimated by the observer, ensuring the exponential stability of the closed-loop system in $H^1$ and $L^2$ norms. Considering now the case of the Neumann boundary measurement (\ref{eq: Neumann output}), the application of the constraints of Theorem~\ref{thm3} are found feasible for $N = 6$ modes estimated by the observer, ensuring the exponential stability of the closed-loop system in $H^1$ norm.

\begin{figure}
     \centering
     	\subfigure[State of the reaction-diffusion system $z(t,x)$]{
		\includegraphics[width=3.5in]{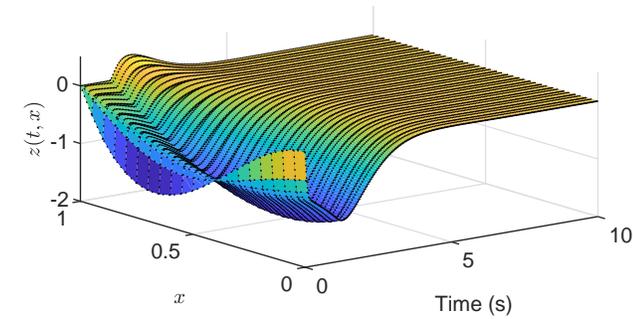}
		}
     	\subfigure[Error of observation $e(t,x) = z(t,x) - \hat{z}(t,x)$]{
		\includegraphics[width=3.5in]{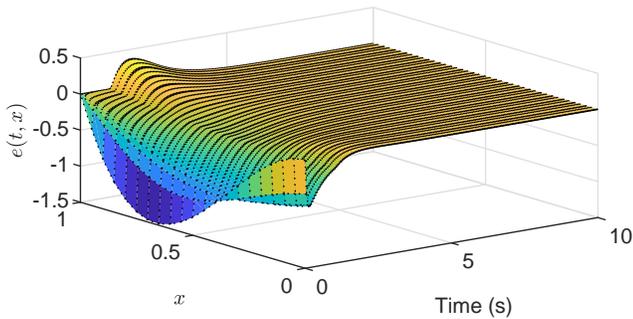}
		}
     	\subfigure[Delayed control input $z(t,1) = u(t-h)$ with delay $h = 1\,\mathrm{s}$]{
		\includegraphics[width=3.5in]{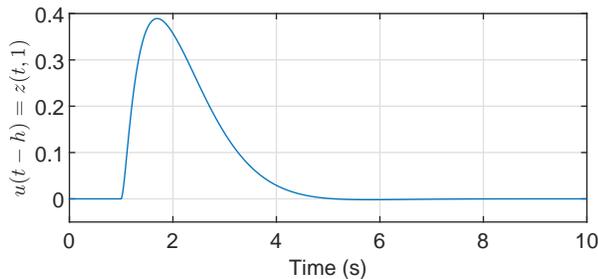}
		\label{fig: sim CL - input}
		}		
     \caption{Time evolution of the closed-loop system for Dirichlet measurement $y_D (t) = z(t,0)$ with delay $h = 1\,\mathrm{s}$}
     \label{fig: sim CL}
\end{figure}

For numerical illustration, we consider the Dirichlet measurement (\ref{eq: Dirichlet output}) along with the initial condition $z_0(x)=10x^2(x-1)$. The evolution of the closed-loop system is depicted in Fig.~\ref{fig: sim CL}. We observe the exponential decay of the both state of the PDE and observation error in spite of the $h = 1\,\mathrm{s}$ input delay. This is compliant with the theoretical prediction of Theorems~\ref{thm1} and~\ref{thm2}.

\section{Conclusion}\label{sec: conclusion}

This paper solved the boundary stabilization problem of general 1-D reaction-diffusion PDEs in the presence of an arbitrarily large input delay. The approach is very general as it covers the cases of Dirichlet/Neumann/Robin boundary control/condition with Dirichlet/Neumann boundary measurement and for system trajectories evaluated in $H^1$ norm (also in $L^2$ norm for Dirichlet measurement). The control strategy couples of finite-dimensional observer that observes a finite number of modes of the PDE and a predictor to compensate the arbitrarily long input delay. Future research directions may be concerned with nonlinear PDEs and non collocated boundary conditions.


\bibliographystyle{plain}        
\bibliography{autosam}           



\appendix
\section{Useful lemma}

The following lemma is borrowed from~\cite[Appendix]{lhachemi2020finite} and is a generalization of a result presented in~\cite{katz2020constructive}.

\begin{lem}\label{lem: useful lemma}
Let $n,m,N \geq 1$, $M_{11} \in \R^{n \times n}$ and $M_{22} \in \R^{m \times m}$ Hurwitz, $M_{12} \in \R^{n \times m}$, $M_{14}^N \in\R^{n \times N}$, $M_{24}^N \in\R^{m \times N}$, $M_{31}^N \in\R^{N \times n}$, $M_{33}^N,M_{44}^N \in \R^{N \times N}$, and
\begin{equation*}
F^N = \begin{bmatrix}
M_{11} & M_{12} & 0 & M_{14}^N \\
0 & M_{22} & 0 & M_{24}^N \\
M_{31}^N & 0 & M_{33}^N & 0 \\
0 & 0 & 0 & M_{44}^N
\end{bmatrix} .
\end{equation*}
We assume that there exist constants $C_0 , \kappa_0 > 0$ such that $\Vert e^{M_{33}^N t} \Vert \leq C_0 e^{-\kappa_0 t}$ and $\Vert e^{M_{44}^N t} \Vert \leq C_0 e^{-\kappa_0 t}$ for all $t \geq 0$ and all $N \geq 1$. Moreover, we assume that there exists a constant $C_1 > 0$ such that $\Vert M_{14}^N \Vert \leq C_1$, $\Vert M_{24}^N \Vert \leq C_1$, and $\Vert M_{31}^N \Vert \leq C_1$ for all $N \geq 1$. Then there exists a constant $C_2 > 0$ such that, for any $N \geq 1$, there exists a symmetric matrix $P^N \in\R^{n+m+2N}$ with $P^N \succ 0$ such that $(F^N)^\top P^N + P^N F^N = - I$ and $\Vert P^N \Vert \leq C_2$.
\end{lem}

\end{document}